\documentclass[a4paper,12pt,final]{amsart}
\usepackage{times,a4wide,mathrsfs, nameref,showkeys,amssymb,amsmath,amsthm,enumerate,xypic,tikzsymbols,dsfont}

\newcommand{\C}{\mathbb{C}}

\newcommand{\QQ}{\mathbb{Q}}
\newcommand{\NN}{\mathbb{N}}
\newcommand{\PP}{\mathbb{P}}

\newcommand{\XX}{\mathcal X}
\newcommand{\YY}{\mathcal Y}

\newcommand{\MM}{\mathcal M}

\newcommand{\rom}{\romannumeral}

\newcommand{\one}{\mathds{1}}

\DeclareMathOperator{\ima}{Im}

\newtheorem{theorem}{Theorem}[section]

\newtheorem{lemma}[theorem]{Lemma}

\newtheorem{corollary}[theorem]{Corollary}
\newtheorem{proposition}[theorem]{Proposition}
\newtheorem{conjecture}[theorem]{Conjecture}
\newtheorem{remark}[theorem]{Remark}
\newtheorem{definition}[theorem]{Definition}
\newtheorem{convention}{Conventions}

\newtheorem{nonumbering}{Theorem}

\newtheorem{nonumberingt}{Acknowledgements}

\begin{document}

\author[Robert Laterveer]
{Robert Laterveer}

\address{Institut de Recherche Math\'ematique Avanc\'ee,
CNRS -- Universit\'e 
de Strasbourg,\
7 Rue Ren\'e Des\-car\-tes, 67084 Strasbourg CEDEX,
FRANCE.}
\email{robert.laterveer@math.unistra.fr}

\title[Questions on the Chow ring of complete intersections]{Questions on the Chow ring of complete intersections}

\begin{abstract} We state several questions, and prove some partial results, about the Chow ring $A^\ast(X)$ of complete intersections in projective space. For one thing, we prove that if $X$ is a general Calabi--Yau hypersurface, the intersection product $A^2(X)\cdot A^i(X)$ is one-dimensional, for any $i>0$. We also show that quintic threefolds have a multiplicative Chow--K\"unneth (MCK) decomposition.
We wonder whether all Calabi--Yau hypersurfaces might have an MCK decomposition, and prove this is the case conditional to a conjecture of Voisin.
\end{abstract}

\thanks{This work is supported by ANR grant ANR-20-CE40-0023.}

\maketitle

\section{Introduction}

Given a complex smooth projective variety $X$, let $A^\ast(X)=\oplus_i A^i(X)$ denote the Chow ring with $\QQ$-coefficients. Even for the simplest varieties, understanding the Chow ring is not so simple.
%there remain many open questions about Chow groups. 
For instance, motivated by the weak Lefschetz theorem in cohomology, Hartshorne has asked the following:

\begin{conjecture}[Hartshorne 1974 \cite{Ha}]\label{ha} Let $X\subset\PP^{n+1}(\C)$ be a smooth hypersurface, and let $h\in A^1(X)$ denote the hyperplane class. Then
  \[ A^i(X)=\QQ[h^i]\ \ \hbox{for\ all}\ i< {n\over 2}\ .\]
  \end{conjecture}
  
 Apart from some easy results when $X$ has small degree (and so is Fano), Conjecture \ref{ha} is completely open for $i\ge 2$, and seems highly challenging. 
  
 Since a direct attack on Conjecture \ref{ha} appears hopeless, let us now investigate some consequences of Conjecture \ref{ha}. As is well-known, the image of intersecting with the hyperplane class on a hypersurface $X$ is one-dimensional, i.e. $h\cdot A^i(X)=\QQ[h^{i+1}]$ (this is just because $h\cdot A^i(X) =\iota^\ast \iota_\ast A^i(X)$ where $\iota\colon X\hookrightarrow \PP^{n+1}(\C)$ denotes the embedding). This observation means that if one believes in Conjecture \ref{ha} one must also believe the following:

\begin{conjecture}\label{ha2} Let $X\subset\PP^{n+1}(\C)$ be a smooth hypersurface. Then
  \[ A^i(X)\cdot A^j(X) =\QQ[h^{i+j}]\ \ \hbox{for\ all}\ i,j >0 \ \hbox{such\ that\ } (i,j)\not=({n\over 2}, {n\over 2})\ .\]
  \end{conjecture}
(It seems likely Conjecture \ref{ha2} holds true more generally for complete intersections $X\subset\PP^{n+r}(\C)$, cf. Remark \ref{cicase} below.)

Restricting attention to hypersurfaces that are Calabi--Yau, there is a remarkable result proven by Voisin:

\begin{theorem}[Voisin \cite{V}]\label{cv} Let $X\subset\PP^{n+1}(\C)$ be a smooth hypersurface of degree $n+2$. Then
  \[ A^i(X)\cdot A^{j}(X) =\QQ[h^{n}]\ \ \ \hbox{for\ all} \   i,j >0 \ \hbox{such\ that\ } i+j=n  \ .\]
  \end{theorem}
  
  Combining Theorem \ref{cv} and Conjecture \ref{ha2}, one obtains the following conjecture about the intersection product on hypersurfaces:
  
  \begin{conjecture}\label{hacy} Let $X\subset\PP^{n+1}(\C)$ be a smooth hypersurface. Assume that either the dimension $n$ is odd, or the degree of $X$ is $n+2$ (i.e. $X$ is Calabi--Yau). Then
  \[ A^i(X)\cdot A^{j}(X) =\QQ[h^{i+j}]\ \ \ \hbox{for\ all} \ i,j>0\ .\]
  \end{conjecture}
  
  By looking into Voisin's proof of Theorem \ref{cv}, we come up with some (very partial) confirmation of Conjecture \ref{hacy}:
  
  \begin{nonumbering}[=Theorem \ref{main}] Let $X\subset\PP^{n+1}(\C)$ be a general hypersurface of degree $n+2$ (i.e. $X$ is Calabi--Yau). Then
  \[ A^2(X)\cdot A^{j}(X) =\QQ[h^{2+j}]\ \ \ \hbox{for\ all} \ j>0\ .\]
  \end{nonumbering}
  
   This settles Conjecture \ref{hacy} for Calabi--Yau hypersurfaces of dimension $\le 6$.  The same result holds for certain Calabi--Yau complete intersections (cf. Theorem \ref{mainc} below). We also prove a result about the ring $B^\ast(X)$ of cycles modulo algebraic equivalence:
  
  \begin{nonumbering}[=Theorem \ref{mainb}]  Let $X\subset\PP^{n+1}(\C)$ be a general hypersurface of degree $n+2$. Then
  \[ B^i(X)\cdot B^{n-i-1}(X) =\QQ[h^{n-1}]\ \ \ \hbox{for\ all}\, 0<i<n-1\ .\]
  \end{nonumbering}
  
  This settles Conjecture \ref{hacy} modulo algebraic equivalence for Calabi--Yau hypersurfaces of dimension $\le 8$.

 Looking at Conjecture \ref{hacy}, it is natural to wonder whether perhaps Calabi--Yau hypersurfaces might have an MCK decomposition, in the sense of Shen--Vial \cite{SV} (roughly speaking, this means that the Chow motive decomposes compatibly with intersection product; cf. Subsection \ref{ss:mck} below for the precise definition). We show this is the case if one assumes a conjecture made by Voisin (cf. Proposition \ref{if} below). We prove an unconditional result in dimension 3:
 
 \begin{nonumbering}[=Theorem \ref{5tic3} and Corollary \ref{cor5tic}] Any smooth quintic threefold $X\subset\PP^4(\C)$ admits an MCK decomposition.
 
 In particular, for any $m\in\NN$ let 
   \[ R^\ast(X^m):=\bigl\langle (p_i)^\ast(h), (p_{jk})^\ast(\Delta_X)\bigr\rangle\ \ \subset\ A^\ast(X^m) \]
   denote the $\QQ$-algebra generated by (pullbacks of) the polarization $h$ and the diagonal $\Delta_X$. Then
   $ R^\ast(X^m)$ injects into cohomology under the cycle class map for all $m\le 205$ (and  $ R^\ast(X^m)$ injects into cohomology for all $m$ if and only if $X$ is Kimura finite-dimensional, in the sense of \cite{Kim}). 
 \end{nonumbering}

 Building on work of Lie Fu \cite{Fu}, we also include a version of Theorem \ref{main} that applies to hypersurfaces that are of general type; this is Theorem \ref{main2} below.

% \begin{theorem}[Fu \cite{Fu}]\label{lf}    Let $X\subset\PP^{n+1}(\C)$ be a general hypersurface of degree $d\ge n+2$. Assume $i_1,\ldots,i_{d-n}$ are strictly positive integers such that $\sum_j i_j=n$. Then
%  \[  A^{i_1}(X)\cdot A^{i_2}(X)\ldots\cdot A^{i_{d-n}}(X)    =\QQ[h^{n}]\ .\] 
%\end{theorem} 
% 

\vskip0.5cm
\begin{convention} In this paper, the word {\sl variety\/} will mean a reduced irreducible scheme of finite type over $\C$. A {\sl subvariety\/} will refer to a (possibly reducible) reduced subscheme which is equidimensional. 

{\bf All Chow groups will be with rational coefficients}: we denote by $A^i(Y)$ the Chow group of codimension $i$ cycles on $Y$ with $\QQ$-coefficients.
The notation $A^i_{hom}(Y)$ (resp. $A^i_{AJ}(X)$) will be used to denote the subgroup of homologically trivial (resp. Abel--Jacobi trivial) cycles.
%For a morphism $f\colon X\to Y$, we write $\Gamma_f\in A^\ast(X\times Y)$ for the graph of $f$.

The contravariant category of Chow motives (i.e., pure motives with respect to rational equivalence and $\QQ$-coefficients as described in \cite{Sc}, \cite{MNP}) will be denoted 
$\MM_{\rm rat}$. 
\end{convention}

\vskip0.5cm

\section{Calabi--Yau hypersurfaces}

\subsection{Intersecting with codimension 2 cycles}

\begin{theorem}\label{main} Let $X\subset\PP^{n+1}(\C)$ be a general hypersurface of degree $n+2$. Then
  \[ A^2(X)\cdot A^i(X) =\QQ[h^{i+2}]\ \ \ \forall i>0\ .\]
  \end{theorem}
  
  \begin{proof} In case $n=4$ this is just Voisin's result (Theorem \ref{cv}). Let us now assume $n\ge 5$.
  
  Let $F:=F(X)$ denote the Fano variety of lines in $X$. By the generality assumption, both $X$ and $F$ are smooth and of the expected dimension, i.e. $\dim F=n-3$.
  Let 
  \[ \Delta_X^{sm}:=\{ (x,x,x)\,\vert\, x\in X\}\ \ \subset\ X\times X\times X\]
  denote the small diagonal. Let
   \begin{equation}\label{gam} \Gamma:= \bigcup_{t\in F} \PP^1_t\times \PP^1_t\times \PP^1_t\ \ \subset\ X\times X\times X\ ,\end{equation}
  where $\PP^1_t\subset X$ denotes the line corresponding to $t\in F$. Let $o_X:={1\over n+2} h^n\in A^n(X)$ denote the ``canonical'' zero-cycle of degree 1, where $h\in A^1(X)$ is the hyperplane section class. For any $1\le i<j\le 3$, define
   \[  \Delta_{ij}:= (p_i\times p_j)^\ast(\Delta_X)\cdot (p_k)^\ast(o_X)\ \ \in\ A^{2n}(X\times X\times X)\ ,\]
   where $k\in \{1,2,3\}\setminus \{i,j\}$, and $p_i$ denotes projection to the $i$th factor.
   
  In the course of the proof of Theorem \ref{cv}, Voisin has obtained \cite[Theorem 3.1]{V} the following equality:
  \begin{equation}\label{sm}    \Delta_X^{sm}={1\over {(n+2)!}} \, \Gamma + \Delta_{12}+\Delta_{13}+\Delta_{23} +P(h_1,h_2,h_3)\ \ \ \hbox{in}\ A^{2n}(X\times X\times X)\ ,\end{equation}
  where $P(h_1,h_2,h_3)$ is a polynomial in the divisors $h_i:=(p_i)^\ast(h)\in A^1(X\times X\times X)$.
  
  Let $a\in A^2(X)$ and $b\in A^i(X)$ (where $i>0$) be any cycles. As $n$ is at least 5, we can write
    \[ a =a_0 + a_2\ \ \hbox{in}\ A^2(X)\ ,\]
    where $\ a_0\in \QQ[h^2]$ and $ a_2\in A^2_{hom}(X)=A^2_{AJ}(X)$.
  Clearly $h\cdot b\in \QQ[h^{i+1}]$, and so it will suffice to prove that $a_2\cdot b$ is in $\QQ[h^{i+2}]$. Considering equality \eqref{sm} as an equality of correspondences from $X\times X$ to $X$, we find that
    \[  a_2\cdot b =  (\Delta_X^{sm})_\ast (a_2\times b) =   \Bigl(  {1\over {(n+2)!}}\,  \Gamma + \Delta_{12}+\Delta_{13}+\Delta_{23} +P(h_1,h_2,h_3)\Bigr){}_\ast (a_2\times b)\ .\]
 The correspondences $\Delta_{12}$ and $P(h_1,h_2,h_3)$ being decomposable, they act as zero on homologically trivial cycles.
 The correspondences $\Delta_{13}$ and $\Delta_{23}$ act as zero on $A^j(X)\otimes A^i(X)$ for all $j,i>0$.
 Thus, the above equality boils down to
  \[     a_2\cdot b =  (\Delta_X^{sm})_\ast (a_2\times b) =     {1\over {(n+2)!}} \, \Gamma_\ast (a_2\cdot b)\ .\]
  Writing $P\in A^{n-1}( F\times X)$ for the (class of the) universal line, we have equality
   \[  \begin{split}\Gamma&= (P\times P\times P)_\ast (\Delta_F^{sm}) \\
                                         &=    P\circ (\Delta_F^{sm}) \circ ({}^t P\times {}^t P)  \ \ \hbox{in}\ A^{2n}(X\times X\times X)\ ,\\
                                         \end{split}\]
                                         where the second equality is an instance of Lieberman's lemma \cite[Lemma 2.1.3]{MNP}. This means that the action of $\Gamma$ on $A^2(X)\otimes A^i(X)$ factors as
              \[  \begin{array}[c]{ccc}
                A^2(X)\otimes A^i(X) & \xrightarrow{\Gamma_\ast} & A^{i+2}(X)\\
                &&\\
               \ \ \ \ \ \ \  \downarrow{\scriptstyle (P^\ast, P^\ast)} &&\ \  \uparrow{\scriptstyle P_\ast}   \\
                &&\\
                A^1(F)\otimes A^{i-1}(F) &\xrightarrow{ (\Delta_F^{sm} )_\ast}& A^{i}(F)\\
                \end{array}\]      
                But since $a_2\in A^2_{AJ}(X)$ we know that $P^\ast(a_2)\in A^1_{AJ}(F)=0$ (here we have used the fact that Abel--Jacobi equivalence is an adequate equivalence relation, hence is compatible with the action of correspondences, cf. \cite[p. 134]{Lie} or \cite[Example 1.11]{Sai}), and so we find
                \[ \Gamma_\ast (a_2\times b)=0\ .\]
                This concludes the proof.                
     \end{proof}

The argument of Theorem \ref{main} can be extended to certain Calabi--Yau complete intersections:
  
  \begin{theorem}\label{mainc} Let $X\subset\PP^{n+r}(\C)$ be a general complete intersection of dimension $n>r$ that is Calabi--Yau, and assume $X$ is not a complete intersection of quadrics. Then
    \[ A^2(X)\cdot A^i(X) =\QQ[h^{i+2}]\ \ \hbox{for\ all}\ \ i>0\ .\]
   \end{theorem} 
   
   \begin{proof} The argument is similar to that of Theorem \ref{main}. 
   In case $n=\dim X=4$ this is a result of Lie Fu \cite[Theorem 0.7]{Fu}. Let us now assume $n\ge 5$.
  
%  Let $F:=F(X)$ denote the Fano variety of lines in $X$. By the generality assumption, both $X$ and $F$ are smooth and of the expected dimension, i.e. $\dim F=n-3$.
%  Let 
%  \[ \Delta_X^{sm}:=\{ (x,x,x)\,\vert\, x\in X\}\ \ \subset\ X\times X\times X\]
%  denote the small diagonal. Let
%   \begin{equation}\label{gam} \Gamma:= \bigcup_{t\in F} \PP^1_t\times \PP^1_t\times \PP^1_t\ \ \subset\ X\times X\times X\ ,\end{equation}
%  where $\PP^1_t\subset X$ denotes the line corresponding to $t\in F$. Let $o_X:={1\over n+2} h^n\in A^n(X)$ denote the ``canonical'' zero-cycle of degree 1, where $h\in A^1(X)$ is the hyperplane section class. For any $1\le i<j\le 3$, define
%   \[  \Delta_{ij}:= (p_i\times p_j)^\ast(\Delta_X)\cdot (p_k)^\ast(o_X)\ \ \in\ A^{2n}(X\times X\times X)\ ,\]
%   where $k\in \{1,2,3\}\setminus \{i,j\}$, and $p_i$ denotes projection to the $i$th factor.
   
  In proving his result, Fu has established \cite[Theorem 1.17]{Fu} the following equality:
  \begin{equation}\label{smf}    \Delta_X^{sm}=\alpha \, \Gamma + (j_{12})_\ast(Z)+  (j_{13})_\ast(Z) +   (j_{23})_\ast(Z)   +P(h_1,h_2,h_3)\ \ \ \hbox{in}\ A^{2n}(X\times X\times X)\ ,\end{equation}
  where $\alpha\in\QQ^\ast$, where $\Gamma$ is defined as in \eqref{gam}, where $j_{12}\colon X\times X\to X^3$ is the diagonal embedding $(x,x^\prime) \mapsto (x,x,x^\prime)$ (and similar definitions for $j_{13}$ and $j_{23}$),
  and $Z:=Q(h_1,h_2)$ is a polynomial in the divisors $h_i:=(p_i)^\ast(h)\in A^1(X\times X)$.
  
  Let $a\in A^2(X)$ and $b\in A^i(X)$ (where $i>0$) be any cycles. As $n$ is at least 5, we can write
    \[ a =a_0 + a_2\ \ \hbox{in}\ A^2(X)\ ,\]
    where $\ a_0\in \QQ[h^2]$ and $ a_2\in A^2_{hom}(X)=A^2_{AJ}(X)$.
  A nice recent result of Mboro \cite[Theorem 1.2]{Mb} ensures that $h\cdot b\in \QQ[h^{i+1}]$ (this uses the assumptions that $n>r$ and that $X$ is not a complete intersection of quadrics), 
  and so it will suffice to prove that $a_2\cdot b$ is in $\QQ[h^{i+2}]$. Considering equality \eqref{smf} as an equality of correspondences from $X\times X$ to $X$, we find that
    \[  \begin{split} a_2\cdot b &=  (\Delta_X^{sm})_\ast (a_2\times b)\\
     &=   \Bigl(  \alpha\,  \Gamma + (j_{12})_\ast(Z)+  (j_{13})_\ast(Z) +   (j_{23})_\ast(Z) +P(h_1,h_2,h_3)\Bigr){}_\ast (a_2\times b)\ .\\
     \end{split}\]
 The correspondences $(j_{12})_\ast(Z)$ and $P(h_1,h_2,h_3)$ being decomposable, they act as zero on homologically trivial cycles.
 The correspondences $(j_{13})_\ast(Z)$ and $(j_{23})_\ast(Z)$ send $A^j(X)\otimes A^i(X)$ to $\QQ[h^{i+j}]$, as can be easily seen (cf. \cite[Proof of Corollary 1.13]{Fu}). 
 Thus, the above equality boils down to
  \[     a_2\cdot b =  (\Delta_X^{sm})_\ast (a_2\times b) =    \alpha \, \Gamma_\ast (a_2\cdot b)\ .\]
  Writing $P\in A^{n-1}( F\times X)$ for the (class of the) universal line, as before we have equality
   \[  \Gamma= 
                                           P\circ (\Delta_F^{sm}) \circ ({}^t P\times {}^t P)  \ \ \hbox{in}\ A^{2n}(X\times X\times X)\ .\]
                                        
            It follows that the action of $\Gamma$ on $A^2(X)\otimes A^i(X)$ factors as
              \[  \begin{array}[c]{ccc}
                A^2(X)\otimes A^i(X) & \xrightarrow{\Gamma_\ast} & A^{i+2}(X)\\
                &&\\
               \ \ \ \ \ \ \  \downarrow{\scriptstyle (P^\ast, P^\ast)} &&\ \  \uparrow{\scriptstyle P_\ast}   \\
                &&\\
                A^1(F)\otimes A^{i-1}(F) &\xrightarrow{ (\Delta_F^{sm} )_\ast}& A^{i}(F)\\
                \end{array}\]      
                But since $a_2\in A^2_{AJ}(X)$ we know that $P^\ast(a_2)\in A^1_{AJ}(F)=0$, and so we conclude that
                \[ \Gamma_\ast (a_2\times b)=0\ .\]
                This ends the proof.                
    \end{proof}
    
%  \begin{remark} One can likewise prove Theorem \ref{mainb} for complete intersections as in Theorem \ref{mainc}; we leave this as an exercice to the reader.
   %  \end{remark}  
     
   \begin{remark}\label{cicase} One might expect that Conjecture \ref{ha2} holds true for all smooth complete intersections in projective space. Indeed, as is well-known Hartshorne's conjecture
   (Conjecture \ref{ha}) for complete intersections would follow from the truth of the Bloch--Beilinson conjectures (cf. \cite[Remark 2.14]{Fu}). Furthermore, one would expect that all smooth complete intersections $X\subset\PP^{\ast}$ satisfy $A^i(X)\cdot h= \QQ[h^{i+1}]$ for all $i$ (indeed, for many complete intersections this is proven in \cite[Theorem 1.2]{Mb} as mentioned above). Combining these two conjectural properties, one obtains Conjecture \ref{ha2} for complete intersections.
     \end{remark}

 \subsection{A result modulo algebraic equivalence} In this subsection, we consider the cycle groups $B^\ast(X):=A^\ast(X)/ A^\ast_{alg}(X)$ (where $A^i_{alg}(X)$ denotes the subgroup of algebraically trivial cycles).
 
 \begin{theorem}\label{mainb}     Let $X\subset\PP^{n+1}(\C)$ be a general hypersurface of degree $n+2$. Then
  \[ B^i(X)\cdot B^{n-i-1}(X) =\QQ[h^{n-1}]\ \ \ \hbox{for\ all}\ \, 0<i<n-1\ .\]
  \end{theorem}
  
  \begin{proof} This is similar to the proof of Theorem \ref{main}.
  The statement being trivially true for $n=4$, we may assume $n\ge 5$.
  Given any $a\in B^i(X)$ and $b\in B^{n-i-1}(X)$, applying once more Voisin's equality \eqref{sm}, their intersection can be expressed as
    \[ a\cdot b =  (\Delta_X^{sm})_\ast (a\times b) =   \Bigl(  {1\over {(n+2)!}}\,  \Gamma + \Delta_{12}+\Delta_{13}+\Delta_{23} +P(h_1,h_2,h_3)\Bigr){}_\ast (a\times b)\ \ \hbox{in}\ B^{n-1}(X)\ .\]
    The correspondences $\Delta_{ij}$ act as zero for dimensional reasons, and it is readily seen that 
      \[ P(h_1,h_2,h_3)_\ast A^\ast(X\times X) =\langle h\rangle\ .\] 
      It follows that
    \[ a\cdot b=    {1\over {(n+2)!}}\,  \Gamma_\ast( a\times b)\ \ \ \hbox{in}\ B^{n-1}(X)\ .\]
  Either $i$ or $n-i-1$ is strictly smaller than $n/2$; without loss of generality, let us assume $i<n/2$. This implies that we can write
   \[  a=a_0 + a_1\ \ \hbox{in}\ B^i(X)  \ ,\]
   where $a_0\in\QQ[h^i]$ and $a_1\in B^i_{hom}(X)$. Clearly $a_0\cdot b\in\QQ[h^{n-1}]$, and so we now restrict attention to the product $a_1\cdot b$.
   As in the proof of Theorem \ref{main} above, the action of $\Gamma$ factors as
     \[  \begin{array}[c]{ccc}
                B^i(X)\otimes B^{n-i-1}(X) & \xrightarrow{\Gamma_\ast} & B^{n-1}(X)\\
                &&\\
               \ \ \ \ \ \ \  \downarrow{\scriptstyle (P^\ast, P^\ast)} &&\ \  \uparrow{\scriptstyle P_\ast}   \\
                &&\\
                B^{i-1}(F)\otimes B^{n-i-2}(F) &\xrightarrow{ (\Delta_F^{sm} )_\ast}& B^{n-3}(F)\\
                \end{array}\]      
                But since $a_1\in B^i_{hom}(X)$ we know that $P^\ast(a_1)\in B^{i-1}_{hom}(F)$ and so 
                \[  P^\ast(a_1)\cdot P^\ast(b)\in B^{n-3}_{hom}(F)=0\ \]
           (indeed, $\dim F=n-3$ and homological and algebraic equivalence coincide for zero-cycles). We conclude that     
                  \[ \Gamma_\ast (a_1\times b)=0\ \ \hbox{in}\ B^{n-1}(X)\ ,\]
   which ends the proof.                
   \end{proof}
   
   \begin{remark} Theorem \ref{mainb} can be extended to complete intersections as in Theorem \ref{mainc}; we leave this as an exercice to the reader.   
    \end{remark}

  \subsection{MCK}\label{ss:mck} 
  
  \begin{definition}[Murre \cite{Mur}] Let $X$ be a smooth projective variety of dimension $n$. We say that $X$ has a {\em CK decomposition\/} if there exists a decomposition of the diagonal
   \[ \Delta_X= \pi^0_X+ \pi^1_X+\cdots +\pi_X^{2n}\ \ \ \hbox{in}\ A^n(X\times X)\ ,\]
  such that the $\pi^i_X$ are mutually orthogonal idempotents and $(\pi_X^i)_\ast H^\ast(X,\QQ)= H^i(X,\QQ)$.
  
  Please note that  ``CK decomposition'' is shorthand for ``Chow--K\"unneth decomposition''.
\end{definition}

\begin{remark} The existence of a CK decomposition for any smooth projective variety is part of Murre's conjectures \cite{Mur}, \cite{J4}. 
\end{remark}

\begin{definition}[Shen--Vial \cite{SV}] Let $X$ be a smooth projective variety of dimension $n$. Let $\Delta_X^{sm}\in A^{2n}(X\times X\times X)$ be the class of the small diagonal
  \[ \Delta_X^{sm}:=\bigl\{ (x,x,x)\ \vert\ x\in X\bigr\}\ \subset\ X\times X\times X\ .\]
  An {\em MCK decomposition\/} is a CK decomposition $\{\pi_X^i\}$ of $X$ that is {\em multiplicative\/}, i.e. it satisfies
  \[ \pi_X^k\circ \Delta_X^{sm}\circ (\pi_X^i\times \pi_X^j)=0\ \ \ \hbox{in}\ A^{2n}(X\times X\times X)\ \ \ \hbox{for\ all\ }i+j\not=k\ .\]
  
Please note that  ``MCK decomposition'' is shorthand for ``multiplicative Chow--K\"unneth decomposition''.
    \end{definition}
    
\begin{remark} Only certain special varieties have an MCK decomposition. For instance, hyperelliptic curves have an MCK, while the general curve of genus $\ge 3$ does {\em not\/} have an MCK. For more on MCK decompositions, cf. \cite{SV}, \cite{FLV2} and the references given there.
\end{remark}    
    
    \subsection{Franchetta property}
    
    \begin{definition}\label{frank} Let $\YY\to B$ be a smooth projective morphism, where $\YY, B$ are smooth quasi-projective varieties. We say that $\YY\to B$ has the {\em Franchetta property in codimension $j$\/} if the following holds: for every $\Gamma\in A^j(\YY)$ such that the restriction $\Gamma\vert_{Y_b}$ is homologically trivial for all $b\in B$, the restriction $\Gamma\vert_{Y_b}$ is zero in $A^j(Y_b)$ for all $b\in B$.
 
 We say that $\YY\to B$ has the {\em Franchetta property\/} if $\YY\to B$ has the Franchetta property in codimension $j$ for all $j$.
 \end{definition}
 
 This property is studied in \cite{FLV}, \cite{FLV3}.
 
 \begin{definition}\label{def:gd} Given a family $\YY\to B$ as above, with $Y:=Y_b$ a fiber, we write
   \[ GDA^j_B(Y):=\ima\Bigl( 
  A^j(\YY)\to A^j(Y)\Bigr) \]
   for the subgroup of {\em generically defined cycles}. 
  In a context where it is clear which family is being referred to, the index $B$ will sometimes be suppressed from the notation.
  \end{definition}
  
  With this notation, the Franchetta property amounts to saying that $GDA^\ast_B(Y)$ injects into cohomology, under the cycle class map. 
  
  \begin{proposition}\label{Fr} Let $B\subset\PP H^0(\PP^{n+1},{\mathcal O}_{\PP^{n+1}}(d))$ be the open subset parameterizing smooth hypersurfaces of degree $d\ge 3$, and let $\XX\to B$ denote the universal family. The families $\XX\to B$ and $\XX\times_B \XX\to B$ have the Franchetta property.
  \end{proposition}
  
  \begin{proof} This is the same argument as \cite[Proposition 5.6]{FLV2}, where this is proven for cubic hypersurfaces.
    \end{proof}

    \subsection{Franchetta and MCK}
    
    As is well-known, complete intersections have a ``standard'' CK decomposition:
    
    \begin{definition}\label{ck} Let $X\subset\PP^{n+r}(\C)$ be a smooth complete intersection of dimension $n$ and degree $d:=\prod_{i=1}^r d_i$. Then
      \[  \pi^i_X:= \begin{cases} {1\over d}\, h^{n-i/2}\times h^{i/2}\ \ \ \ \hbox{if}\ i\not=n\ \hbox{and}\ i\ \hbox{is\ even}\ ,\\
                                            0\ \ \ \ \hbox{if}\ i\not=n\ \hbox{and}\ i\ \hbox{is\ odd}\ ,\\
                                        \Delta_X-\sum_{i\not= n} \pi^i_X\ \ \ \ \hbox{if}\ i=n\ \\
                                        \end{cases}\]
               defines a CK decomposition.                         
      \end{definition}
      
    \begin{remark} In the set-up of Definition \ref{ck}, in case $n$ is even one can further decompose
      \[ \pi^n_X= \pi^{n,alg}_X + \pi^{n,prim}_X\ ,\]
      where $\pi^{n,alg}_X:=  {1\over d}\, h^{n/2}\times h^{n/2}$.
      \end{remark}

    \begin{corollary}\label{only} Let $X\subset\PP^{n+1}(\C)$ be a smooth hypersurface of degree at least 3, and $\{\pi^i_X\}$ as in Definition \ref{ck}. Assume that
      \begin{equation}\label{mot} \pi^{n,prim}_X\circ \Delta_X^{sm}\circ (  \pi^{n,prim}_X\times \pi^{n,prim}_X) =0\ \ \hbox{in}\ A^{2n}(X\times X\times X)\ .\end{equation}
     Then $\{\pi^i_X\}$ is an MCK decomposition     
    \end{corollary}
    
    \begin{proof} The point is that (by the very construction of the projectors $\pi^i_X$) one has isomorphisms of Chow motives
      \begin{equation}\label{one} (X,\pi_X^{2i},0)=\one(-i)\ \ \forall\ i\not=n\ ,\ \ (X,\pi_X^{n,alg},0)=\one(-{n\over 2})\ \ \ \hbox{in}\ \MM_{\rm rat}\ .\end{equation}
      
      Let $\pi^i_X,\pi^j_X,\pi^k_X$ be three projectors such that $i+j\not=k$, so that
      \[    \pi_X^k\circ \Delta_X^{sm}\circ (\pi_X^i\times \pi_X^j)\ \in\ A^{2n}_{hom}(X\times X\times X)    \ .\]
     In view of \eqref{mot}, we may assume that at least one of  $\pi^i_X,\pi^j_X,\pi^k_X$ is different from $\pi^{n,prim}_X$. 
     Then we have that
        \[ \begin{split} \pi_{X}^k\circ \Delta_{X}^{sm}\circ (\pi_{X}^i\times \pi_{X}^j)&= ({}^t     \pi_{X}^i\times {}^t \pi_{X}^j\times \pi^k_{X})_\ast \Delta^{sm}_{X}\\
                                                           & = ({}   \pi_{X}^{2n-i}\times {}\pi_{X}^{2n-j}\times \pi^k_{X})_\ast \Delta^{sm}_{X}\\  
                                                           & \hookrightarrow \bigoplus A^\ast(X\times X)\ .\\
                                                           \end{split}\]
                                                           Here the first equality is an application of Lieberman's lemma \cite[Lemma 2.1.3]{MNP}, the second equality is by self-duality of the $\{\pi^\ast_X\}$, and the inclusion follows from property \eqref{one}. The resulting cycle in $\bigoplus A^\ast(X\times X)$ is generically defined (since the $\pi^\ast_{X}$ and $\Delta_{X}^{sm}$ are generically defined) and homologically trivial (since $i+j\not=k$). By the Franchetta property for $X\times X$ (Proposition \ref{Fr}), the resulting cycle in $\bigoplus A^\ast(X\times X)$ is rationally trivial, and so        
                                                          \[ \pi_{X}^k\circ \Delta_{X}^{sm}\circ (\pi_{X}^i\times \pi_{X}^j)=0\ \ \ \hbox{in}\ A^{2n}(X\times X\times X)\ ,\]
                                                 as desired. This proves the corollary.
                                                 \end{proof}

     %Then \eqref{mot} implies that 
   %  \[   \pi^k_X\circ [ -]     \circ ( \pi_X^i\times \pi_X^j)  \colon\ \ A^{\ast}_{}(X\times X\times X)\ \to\     A^{\ast}_{}(X\times X\times X) \]
  %   factors over $A^\ast(X\times X)$, and in particular the restriction
    %  \[  \pi^k_X\circ [ -]     \circ ( \pi_X^i\times \pi_X^j)  \colon\ \ GDA^{\ast}_{hom}(X\times X\times X)\ \to\     GDA^{\ast}_{hom}(X\times X\times X) \]
    % factors over $GDA^\ast_{hom}(X\times X)$. 
  %  By Lieberman's lemma, we have equality
   %   \[   \pi_X^k\circ \Delta_X^{sm}\circ (\pi_X^i\times \pi_X^j)= (\pi_X )_\ast \Delta_X^{sm}    \ \ \hbox{in}\ A^{2n}_{hom}(X\times X\times X)    \ .\]

%          
%    But $GDA^\ast_{hom}(X\times X)=0$ by Proposition \ref{Fr}, and so 
%     \[  \pi^k_X\circ GDA^\ast_{hom}(X\times X\times X)    \circ  ( \pi_X^i\times \pi_X^j) =0\ .\]
%     Applying this to   $\pi_X^k\circ \Delta_X^{sm}\circ (\pi_X^i\times \pi_X^j) \in A^{2n}_{hom}(X\times X\times X)    $, we find that
%     \[    \pi_X^k\circ \Delta_X^{sm}\circ (\pi_X^i\times \pi_X^j)=0\ \hbox{in}\ A^{2n}_{}(X\times X\times X)   \ ,\]
%     proving the corollary.                
%    \end{proof}
  
    \subsection{A conditional result}
  
  In her paper proving Theorem \ref{cv}, Voisin has stated the following conjecture:
  
  \begin{conjecture}[Voisin \cite{V}]\label{conj} Let $X\subset\PP^{n+1}(\C)$ be a general hypersurface of degree $n+2$, and let $\Gamma\in A^{2n}(X\times X\times X)$ be the cycle defined in \eqref{gam}. Then
    \[ \Gamma\in \ima\bigl( A^\ast(\PP^{n+1}\times\PP^{n+1}\times\PP^{n+1})\to A^\ast(X\times X\times X)\bigr)\ .\]
     \end{conjecture}  
     
    (This is \cite[Conjecture 3.5]{V}.)
  
  Voisin's conjecture has strong consequences:
    
   \begin{proposition}\label{if} Assume Conjecture \ref{conj}. Let  $X\subset\PP^{n+1}(\C)$ be a smooth hypersurface of degree $n+2$.
   
   \noindent
   (\rom1) Conjecture \ref{hacy} is true for $X$, i.e. 
       \[ A^i(X)\cdot A^j(X) =\QQ[h^{i+j}]\ \ \hbox{for\ all}\ \ i,j>0\ .\]
       
   \noindent
   (\rom2) $X$ has an MCK decomposition.
   \end{proposition}    
   
   \begin{proof} Conjecture \ref{conj} together with the equality \eqref{sm} imply that for a general hypersurface $X$ there is equality
   \begin{equation}\label{sm1} 
     \Delta_X^{sm}= \Delta_{12}+\Delta_{13}+\Delta_{23} +P(h_1,h_2,h_3)\ \ \ \hbox{in}\ A^{2n}(X\times X\times X)\ ,\end{equation}
   where  $P(h_1,h_2,h_3)$ is a new polynomial in the divisor classes $h_i$. In view of the spread lemma \cite[Lemma 3.2]{Vo}, equality \eqref{sm1} then holds for {\em all\/} smooth Calabi--Yau hypersurfaces $X$. This immediately implies (\rom1) (cf. the proof of Theorem \ref{main} above).

 To prove (\rom2), in view of Corollary \ref{only}, we just need to check the vanishing
   \begin{equation}\label{??}  \pi^{n,prim}_X\circ \Delta_X^{sm}\circ (   \pi^{n,prim}_X\times \pi^{n,prim}_X){\stackrel{}{=}}\, 0\ .\end{equation}
 In view of equality \eqref{sm1}, we can write
  \[     \pi^{n,prim}_X\circ \Delta_X^{sm}\circ (   \pi^{n,prim}_X\times \pi^{n,prim}_X)=     \pi^{n,prim}_X\circ \bigl( \sum \Delta_{ij}+P(h_1,h_2,h_3)\bigr)   \circ (   \pi^{n,prim}_X\times \pi^{n,prim}_X)  \ .\]
  
  As for the $\Delta_{ij}$, it follows from their construction that there is equality
    \[    \Delta_{12}= \pi^{2n}_X\circ \Delta_{12}\ ,\ \ \Delta_{13}= \Delta_{13}\circ (\Delta_X\times \pi^0_X)\ ,\ \ \Delta_{23}= \Delta_{23}\circ (\pi^0_X\times\Delta_X)\ .\]
    The projectors $\pi^j_X$ being orthogonal, we thus find that
    \[    \pi^{n,prim}_X\circ ( \sum \Delta_{ij})   \circ (   \pi^{n,prim}_X\times \pi^{n,prim}_X)    =0\ .\]
    Likewise, any monomial in the $h_i$ satisfies
    \[    (h_1^{i} h_2^j h_3^{2n-i-j})=  \pi_X^{i+j-n,alg}  \circ (h_1^{i} h_2^j h_3^{2n-i-j})\circ (\pi_X^{2n-i,alg}\times \pi_X^{2n-j,alg})    \]
  (NB: the superscript ``alg'' is only operative for $\pi_X^n$, $n$ even, i.e. we use the convention $\pi_X^{i,alg}=\pi_X^i$ for $i\not=n$ ).  
  Again by orthogonality of the projectors, we thus find that
    \[      \pi^{n,prim}_X\circ \bigl( P(h_1,h_2,h_3)\bigr)   \circ (   \pi^{n,prim}_X\times \pi^{n,prim}_X) =0 \ .\]
    This proves the required vanishing \eqref{??}, and ends the proof of the proposition.
     \end{proof}

   \subsection{Quintic threefolds}
   
   In dimension 3, we can prove an unconditional result:
   
   \begin{theorem}\label{5tic3} Any smooth quintic threefold $X\subset\PP^4(\C)$ admits an MCK decomposition.
   \end{theorem}
   
   \begin{proof} The point is that for the general quintic threefold, the Fano variety of lines $F=F(X)$ is zero-dimensional (more precisely, $F$ consists of 2875 reduced points \cite{Sch1}, \cite{Sch2}, \cite[Corollary 6.35]{EH}).
   
   To construct an MCK decomposition for the general quintic threefold, in view of Corollary \ref{only}, we just need to check the vanishing
   \begin{equation}\label{?}  \pi^{3}_X\circ \Delta_X^{sm}\circ (   \pi^{3}_X\times \pi^{3}_X){\stackrel{}{=}}\, 0\ \ \hbox{in}\ A^6(X\times X\times X)\ .\end{equation}
 In view of equality \eqref{sm}, we can write
  \[     \pi^{3}_X\circ \Delta_X^{sm}\circ (   \pi^{3}_X\times \pi^{3}_X)=     \pi^{3}_X\circ \bigl({1\over {5!}} \, \Gamma +\sum \Delta_{ij}+P(h_1,h_2,h_3)\bigr)   \circ (   \pi^{3}_X\times \pi^{3}_X)  \ .\]
  The summands involving $\Delta_{ij}$ and $P(h_1,h_2,h_3)$ vanish for general reasons (cf. the proof of Proposition \ref{if}(\rom2) above), and so it only remains to analyze the summand involving $\Gamma$. As before, we can write
    \[      \Gamma=     P\circ (\Delta_F^{sm}) \circ ({}^t P\times {}^t P)  \ \ \hbox{in}\ A^{6}(X\times X\times X)   \ ,\]
    where $P\in A^2 (F\times X)$ is the (class of the) universal line. In particular, we have
    \[    \pi^3_X\circ \Gamma\circ (   \pi^{3}_X\times \pi^{3}_X)=   \pi^3_X\circ  P\circ (\Delta_F^{sm}) \circ ({}^t P\times {}^t P) \circ (   \pi^{3}_X\times \pi^{3}_X) \ \ \hbox{in}\ A^{6}(X\times X\times X)   \ .\]  
  But now, Lemma \ref{P} below (combined with the orthogonality $\pi^3_X\circ \pi^4_X=0$) implies the vanishing
   \[    \pi^3_X\circ \Gamma\circ (   \pi^{3}_X\times \pi^{3}_X)=  0\ \     \hbox{in}\ A^{6}(X\times X\times X)   \ .\]  
   This shows \eqref{?} and gives an MCK decomposition for the general quintic threefold. 
   
   To extend to {\em all\/} smooth quintic threefolds, one observes that all terms in \eqref{?} are generically defined, and so the spread lemma \cite[Lemma 3.2]{Vo} implies the vanishing \eqref{?} for all smooth quintic threefolds.
   
   \begin{lemma}\label{P} Let $X$ be a general quintic threefold, and  $P\in A^2 (F\times X)$ the universal line. 
   Then
     \[   P= \pi^4_X\circ P\ \ \ \hbox{in}\ A^2(F\times X)\ .\]
     \end{lemma}
     
     To prove the lemma, we observe that the equality is true modulo homological equivalence (indeed, the correspondence $P$ sends $H^\ast(F,\QQ)=H^0(F,\QQ)$ to $H^4(X,\QQ)$). All terms being generically defined (and remembering that $A^2(F\times X) =A^2(X)^{\oplus 2875}$), the lemma then follows from the Franchetta property for $X$ (Proposition \ref{Fr}).
        \end{proof}
        
  \begin{remark} Unfortunately, the argument proving Theorem \ref{5tic3} breaks down for Calabi--Yau hypersurfaces of dimension greater than $3$. The reason is that in Lemma \ref{P}, the argument hinges on the Franchetta property for $F\times X$. This is easy when the dimension of $F$ is zero, but becomes problematic as soon as the dimension of $F$ is greater than $0$.
  \end{remark}

   \begin{corollary}\label{cor5tic} Let $X$ be a smooth quintic threefold. For any $m\in\NN$ let 
   \[ R^\ast(X^m):=\bigl\langle (p_i)^\ast(h), (p_{jk})^\ast(\Delta_X)\bigr\rangle\ \ \subset\ A^\ast(X^m) \]
   denote the $\QQ$-algebra generated by (pullbacks of) the polarization $h$ and the diagonal $\Delta_X$. Then
   $ R^\ast(X^m)$ injects into cohomology under the cycle class map for all $m\le 205$. Moreover,  $ R^\ast(X^m)$ injects into cohomology for all $m$ if and only if $X$ is Kimura finite-dimensional, in the sense of \cite{Kim}.    
   \end{corollary}
   
   \begin{proof} This is an application of \cite[Proposition 2.11]{FLV3}, using the fact that $X$ has an MCK decomposition and that $\dim H^3(X,\QQ)=204$.
     \end{proof}
     
     \begin{remark} Conjecturally, all smooth projective varieties are Kimura finite-dimensional \cite{Kim}. An example which is known to be Kimura finite-dimensional is the Fermat quintic threefold 
       \[ x_0^5+\dots +x_4^5=0 \ .\] 
 (This can be proven using Shioda's inductive structure of Fermat hypersurfaces, cf. \cite{SI}.)
     \end{remark}
     
  \begin{remark} As explained in \cite[Section 2.3]{FLV3}, the statement of Corollary \ref{cor5tic} is inspired by results on the so-called ``tautological ring'' of hyperelliptic curves \cite{Tav0}, \cite{Tav} and of K3 surfaces \cite{V17}, \cite{Yin}.
   \end{remark}

\section{General type hypersurfaces} In this section, we consider hypersurfaces of degree $d$ strictly larger than $n+2$; these hypersurfaces are of general type. Hartshorne's conjecture (Conjecture \ref{ha}) has the following consequence:

\begin{conjecture}\label{3} Let $X\subset\PP^{n+1}(\C)$ be a smooth hypersurface. Then
  \[ A^{i_1}(X)\cdot A^{i_2}(X)\cdot A^{i_3}(X) =\QQ[h^{i_1+i_2+i_3}]\ \ \hbox{for\ all}\ \ i_1,i_2,i_3>0\ .\]
  \end{conjecture}

\subsection{Lie Fu's theorem} Inspired by Voisin's result for Calabi--Yau hypersurfaces (Theorem \ref{cv}), Lie Fu has proven a nice result about zero-cycles that are intersections on general type hypersurfaces, providing partial confirmation to Conjecture \ref{3}:

 \begin{theorem}[Fu \cite{Fu}]\label{lf}    Let $X\subset\PP^{n+1}(\C)$ be a general hypersurface of degree $d\ge n+2$. Assume $i_1,\ldots,i_{d-n}$ are strictly positive integers such that $\sum_j i_j=n$. Then
  \[  A^{i_1}(X)\cdot A^{i_2}(X)\ldots\cdot A^{i_{d-n}}(X)    =\QQ[h^{n}]\ .\] 
\end{theorem} 
 
\subsection{Main result} The main result in this section is that in the setting of Theorem \ref{lf}, one can obtain a stronger result if one of the codimensions $i_j$ is equal to 2:

\begin{theorem}\label{main2} Let $X\subset\PP^{n+1}(\C)$ be a general hypersurface of degree $d\ge n+2$. Then
  \[ A^2(X)\cdot A^{i_1}(X)\cdot A^{i_2}(X)\cdot \ldots\cdot A^{i_{d-n-1}}(X)    =\QQ[h^{2+\sum_j i_j}]\ \ \ \hbox{for\ all}\ \, i_1,\ldots , i_{d-n-1}>0\ .\]
\end{theorem}

\begin{proof} Let $k:=d+1-n$ (and so $k\ge 3$), and let
  \[  \delta_X:=\{(x,x,\ldots,x)\vert x\in X\}\ \ \subset \ X^k\]
  denote the smallest diagonal. By the generality assumption, the Fano variety $F:=F(X)$ of lines in $X$ is smooth of dimension $n-k$.
In the course of proving Theorem \ref{lf}, Fu has obtained \cite[Theorem 2.12]{Fu} the following dichotomy for $\delta_X$:
either
  \begin{equation}\label{either} \delta_X= {(-1)^{k-1}\over d!}\, \Gamma +\sum_{i=1}^k D_i + \sum_j \lambda_j \sum_{\vert I\vert=j} D_I + P(h_1,\ldots,h_k)\ \ \hbox{in}\ A^\ast(X^k)\ ,\end{equation}
  or there exists $\ell<k$ such that
  \begin{equation}\label{or} \delta_X= \sum_{i=1}^k D_i + \sum_j \lambda_j \sum_{\vert I\vert=j} D_I + P(h_1,\ldots,h_k)\ \ \hbox{in}\ A^\ast(X^\ell)\ .\end{equation}
Here, $\Gamma$ is defined as
  \[ \Gamma:= \bigcup_{t\in F(X)} \PP^1_t\times\cdots\times\PP^1_t\ \ \subset\ X^k\ ,\]
  the cycle $D_i$ is defined as $(p_i)^\ast(o_X)\cdot \Delta_{i^c}$, where $  \Delta_{i^c}$ is the diagonal of the complementary set $\{1,\ldots,k\}\setminus i$, and similarly, for any subset $I\subset\{1,\ldots,k\}$, the cycle
$D_I$ is defined as $\prod_{i\in I}(p_i)^\ast(o_X)\cdot \Delta_{I^c}$.

Once again, we consider the equality \eqref{either} (resp. \eqref{or}) as an equality of correspondences from $X^{k-1}$ (resp. $X^{\ell-1}$) to $X$.
Given $d-n$ cycles of positive codimension $a_0,\ldots,a_{d-n-1}\in A^\ast(X)$, it is readily seen that 
  \[  (D_i)_\ast(a_0\times\cdots\times a_{d-n-1})=  (D_I)_\ast(a_0\times\cdots\times a_{d-n-1})  =0\ \ \hbox{in}\ A^\ast(X)\ ,\]
  while
  \[  \bigl( P(h_1,\ldots,h_k)\bigr){}_\ast(a_0\times\cdots\times a_{d-n-1})\ \ \in\ \QQ[h^{\ast}]\ .\]
 Hence, to prove Theorem \ref{main2} one only needs to worry about the action of $\Gamma$. Since (by Lieberman's lemma) $\Gamma$ can be written as
   \[  \Gamma=   P\circ (\delta_F^{}) \circ ({}^t P\times\cdots\times  {}^t P)  \ \ \hbox{in}\ A^{\ast}(X^k)\ ,
                                       \]
             the action of $\Gamma$ factors as
              \[  \begin{array}[c]{ccc}
                A^2(X)\otimes A^{i_1}(X) \otimes\cdots\otimes A^{i_{k-2}}(X)& \xrightarrow{\Gamma_\ast} & A^{2+\sum_j i_j}(X)\\
                &&\\
               \ \ \ \ \ \ \  \downarrow{\scriptstyle (P^\ast, \ldots, P^\ast)} &&\ \  \uparrow{\scriptstyle P_\ast}   \\
                &&\\
                A^1(F)\otimes A^{i_1-1}(F) \otimes\cdots\otimes A^{i_{k-2}-1}(F)&\xrightarrow{ (\delta_F^{} )_\ast}& A^{\sum_j i_j -k+3}(F)\\
                \end{array}\]      
  The result now follows from the fact that $A^1_{AJ}(F)=0$ (cf. the proof of Theorem \ref{main}).
   \end{proof}

  \begin{remark} Similarly to Theorem \ref{mainb}, one can prove a result for one-cycles modulo algebraic equivalence for general type hypersurfaces; we leave this as an exercice for the diligent reader.
   \end{remark}

% 
% \vskip0.5cm
%\begin{nonumberingda} The author states that this is not applicable, as there are no associated data.
%\end{nonumberingda}
%
%\vskip0.5cm
%\begin{nonumberingcon} The author states that there is no conflict of interest.
%\end{nonumberingcon}

\vskip0.5cm

 \begin{nonumberingt} Thanks to the highly efficient coffee machine at the Schiltigheim Research Center in Advanced Mathematics.
\end{nonumberingt}

\end{document}